\documentclass[12pt]{article}
\usepackage{amssymb}
\usepackage{amsfonts}
\usepackage{amsmath}

\input{tcilatex}
\begin{document}

\begin{center}
\textbf{A generalization of Istr\u{a}\c{t}escu's fixed point theorem}

\textbf{for convex contractions}

\bigskip

Radu MICULESCU and Alexandru MIHAIL

\bigskip
\end{center}

\textbf{Abstract}. {\small In this paper we prove a generalization of Istr%
\u{a}\c{t}escu's theorem for convex contractions. More precisely, we
introduce the concept of iterated function system consisting of convex
contractions and prove the existence and uniqueness of the attractor of such
a system. In addition we study the properties of the canonical projection
from the code space into the attractor of an iterated function system
consisting of convex contractions.}

\bigskip


{\small iterated function systems consisting of convex contractions}

\bigskip

\textbf{1. Introduction}

\bigskip

Banach-Caccioppoli-Picard contraction principle, which is an extremely
useful tool in nonlinear analysis, says that any contraction $%
f:(X,d)\rightarrow (X,d)$, where $(X,d)$ is a complete metric space, has a
unique fixed point $x^{\ast }$ and $\underset{n\rightarrow \infty }{\lim }%
f^{[n]}(x)=x^{\ast }$ for every $x\in X$. Besides its great features (the
uniqueness of the fixed point and the possibility to approximate it by the
means of Picard iteration) there exists a drawback of this result, namely
that the contraction condition is too strong.

The natural question if there exist contraction-type conditions that do no
imply the contraction condition and for which the existence and uniqueness
of the fixed point are assured was answered, among others, by V. Istr\u{a}%
\c{t}escu who introduced and studied the convex contraction condition (see
[5], [6] and [7]). More precisely, a continuous function $f:(X,d)\rightarrow
(X,d)$, where $(X,d)$ is a complete metric space, is called convex
contraction if there exist $a,b\in (0,1)$ such that $a+b<1$ and $%
d(f^{[2]}(x),f^{[2]}(y))\leq ad(f(x)),f(y))+bd(x,y)$ for every $x,y\in X$.
Istr\u{a}\c{t}escu proved that any convex contraction has a unique fixed
point $x^{\ast }\in X$ (and $\underset{n\rightarrow \infty }{\lim }%
f^{[n]}(x)=x^{\ast }$ for every $x\in X$) and provided an example of convex
contraction which is not contraction. V. Ghorbanian, S. Rezapour and N.
Shahzad [8] generalized Istr\u{a}\c{t}escu's results to complete ordered
metric spaces. M. A. Miandaragh, M. Postolache and S. Rezapour [16]
introduced the concept of generalized convex contraction and proved some
theorems about approximate fixed points of these contractions. Extending
these results, A. Latif, W. Sintunavarat and A. Ninsri [12] introduced a new
concept called partial generalized convex contraction and established some
approximate fixed point results for such mappings in $\alpha $-complete
metric spaces. For more results along these lines of generalization one can
also see [10].

Let us recall that an iterated function system on a complete metric space%
\textit{\ }$(X,d)$, denoted by $\mathcal{S}=(X,(f_{k})_{k\in \{1,2,...,n\}})$%
, consists of a finite family of contractions $(f_{k})_{k\in \{1,2,...,n\}}$%
, where $f_{k}:X\rightarrow X$. The function $F_{\mathcal{S}}:\mathcal{K}%
(X)\rightarrow \mathcal{K}(X)$ defined by $F_{\mathcal{S}}(C)=\underset{k=1}{%
\overset{n}{\cup }}f_{k}(C)$, for all $C\in \mathcal{K}(X)$ -the set of
non-empty compact subsets of $X$-, which is called the set function
associated to $\mathcal{S}$, turns out to be a contraction (with respect to
the Hausdorff-Pompeiu distance) and its unique fixed point, denoted by%
\textit{\ }$A_{\mathcal{S}}$\textit{, }is called the attractor of the system%
\textit{\ }$\mathcal{S}$. As iterated function systems represent one of the
main tools to generate fractals, the extending problem of the notion of
iterated function system was treated by several authors. Let us mentions
some contributions along these lines of research. Given a complete metric
space $(X,d)$ and a finite family of functions $f_{1},f_{2},...,f_{n}:X%
\rightarrow X$, L. M\'{a}t\'{e} [15] proved the existence of a unique $A\in 
\mathcal{K}(X)$ such that $A=\underset{i=1}{\overset{n}{\cup }}f_{i}(A)$
under weaker contractivity conditions (for example $d(f_{i}(x),f_{i}(y))\leq
\varphi (d(x,y))$, where $\varphi :[0,\infty )\rightarrow \lbrack 0,\infty )$
is an upper continuous non-decreasing function with the property that $%
\varphi (t)<t$ for each $t>0$). K. Le\'{s}niak [13] presented a multivalued
approach of infinite iterated function systems. A. Petru\c{s}el [21] proved
that each finite family of single-valued and multi-valued operators
satisfying some Meir-Keeler type conditions has a self-similar set (see also
[4]). Let $(X,d)$ be a metric space and $f_{1},f_{2},...,f_{n}:X\rightarrow
P_{cl}(X)$ be set-valued mappings on $X$, where $P_{cl}(X)$ designates the
family of all nonempty closed subsets of $X$. The system $%
F=(f_{1},f_{2},...,f_{n})$ is called an iterated multifunction system and
the operator $\overset{\symbol{94}}{F}:P_{cl}(X)\rightarrow P_{cl}(X)$ given
by $\overset{\symbol{94}}{F}(Y)=\overline{\underset{i=1}{\overset{n}{\cup }}%
f_{i}(Y)}$, where $f_{i}(Y)=\underset{x\in Y}{\cup }f_{i}(x)$ for each $i\in
\{1,2,...,n\}$, is called the Barnsley-Hutchinson operator generated by $F$.
A fixed point of this operator is called a multivalued large fractal. C.
Chifu and A. Petru\c{s}el [3] obtained existence and uniqueness results for
multivalued large fractals (see also [20]). G. Gw\'{o}\'{z}d\'{z}-\L ukowska
and J. Jachymski [9] developed the Hutchinson-Barnsley theory for finite
families of mappings on a metric space endowed with a directed graph. E.
Llorens-Fuster, A. Petru\c{s}el and J.-C. Yao [14] gave existence and
uniqueness results for self-similar sets of a mixed iterated function
system. M. Boriceanu, M. Bota and A. Petru\c{s}el [2] extended the
Hutchinson-Barnsley theory to the case of set-valued mappings on a $b$%
-metric space. For other related results see [1], [11], [17], [19], [22],
[24], [25] and [26].

In this paper we introduce the concept of iterated function system
consisting of convex contractions and prove the existence and uniqueness of
the attractor of such a system obtaining in this way a generalization of Istr%
\u{a}\c{t}escu's convex contractions fixed point theorem (see Theorem 3.2).
Moreover we study the properties of the canonical projection from the code
space into the attractor of an iterated function system consisting of convex
contractions (see Theorem 3.6).

\bigskip

\textbf{2. Preliminaries}

\bigskip

Given a function $f:X\rightarrow X$, by $f^{[n]}$ we mean the composition of 
$f$\ by itself $n$\ times.

Given a set $X$ and a family of functions $(f_{i})_{i\in I}$, where $%
f_{i}:X\rightarrow X$, by $f_{\alpha _{1}\alpha _{2}....\alpha _{n}}$ we
mean $f_{\alpha _{1}}\circ f_{\alpha _{2}}\circ ...\circ f_{\alpha _{n}}$
and by $Y_{\alpha _{1}\alpha _{2}....\alpha _{n}}$ we understand $f_{\alpha
_{1}\alpha _{2}....\alpha _{n}}(Y)$, where $Y\subseteq X$ and $\alpha
_{1},\alpha _{2},....,\alpha _{n}\in I$.

Given a set $X$, by $\mathcal{P}^{\ast }(X)$ we denote the family of all
nonempty subsets of $X$. For a metric space $(X,d)$, by $\mathcal{K}(X)$ we
denote the set of non-empty compact subsets of $X$.

Given two sets $A$ and $B$, by $B^{A}$ we mean the set of functions from $A$
to $B$.

Given a set $I$, $\Lambda (I)$ denotes $I^{\mathbb{N}^{\ast }}$ and $\Lambda
_{n}(I)$ denotes $I^{\{1,2,...,n\}}$. Hence the elements of $\Lambda (I)$
can be written as infinite words $\alpha =\alpha _{1}\alpha _{2}\alpha
_{3}...$ and the elements of $\Lambda _{n}(I)$ as finite words $\alpha
=\alpha _{1}\alpha _{2}....\alpha _{n}$. By $\Lambda ^{\ast }(I)$ we denote
the set of all finite words, namely $\Lambda ^{\ast }(I)=\underset{n\in 
\mathbb{N}^{\ast }}{\cup }\Lambda _{n}(I)\cup \{\lambda \}$, where $\lambda $
is the empty word. $\Lambda (I)$ can be seen as a metric space with the
distance $d_{\Lambda }$ defined by $d_{\Lambda }(\alpha ,\beta )=\frac{1}{%
2^{n}}$ where $n$ is the natural number having the property that $\alpha
_{k}=\beta _{k}$ for $k<n$ and $\alpha _{n}\neq \beta _{n}$ if $\alpha
=\alpha _{1}\alpha _{2}\alpha _{3}...\alpha _{n}\alpha _{n+1}...\neq \beta
=\beta _{1}\beta _{2}\beta _{3}...\beta _{n}\beta _{n+1}...$ and $d_{\Lambda
}(\alpha ,\alpha )=0$. By $\alpha \beta $ we understand the concatenation of
the words $\alpha \in \Lambda ^{\ast }$ and $\beta \in \Lambda \cup \Lambda
^{\ast }$. For $\alpha \in \Lambda \cup \Lambda _{n}$ and $m\leq n$, $%
[\alpha ]_{m}\overset{def}{=}\alpha _{1}\alpha _{2}....\alpha _{m}$. For $%
i\in I$, let us consider the function $F_{i}:\Lambda (I)\rightarrow \Lambda
(I)$\ given by $F_{i}(\alpha )=i\alpha $ for all $\alpha \in \Lambda (I)$.

\newpage

\textbf{Definition 2.1.} \textit{For a metric space }$(X,d)$\textit{, we
consider on} $\mathcal{P}^{\ast }(X)$ \textit{the generalized
Hausdorff-Pompeiu pseudometric} $h:\mathcal{P}^{\ast }(X)\times \mathcal{P}%
^{\ast }(X)\rightarrow \lbrack 0,+\infty ]$ \textit{defined by}

\begin{equation*}
h(A,B)=\max (d(A,B),d(B,A))=\inf \{\eta \in \lbrack 0,\infty ]\mid
A\subseteq N_{\eta }(B)\text{ \textit{and} }B\subseteq N_{\eta }(A)\}
\end{equation*}%
\textit{where }%
\begin{equation*}
d(A,B)=\underset{x\in A}{\sup }d(x,B)=\underset{x\in A}{\sup }(\underset{%
y\in B}{\inf }d(x,y))
\end{equation*}
\textit{\ and }%
\begin{equation*}
N_{\eta }(A)=\{x\in X\mid \text{\textit{there exists}}\mathit{\ }y\in X%
\mathit{\ }\text{\textit{such that }}d(x,y)<\eta \}\text{,}
\end{equation*}%
\textit{for every} $A,B\in \mathcal{P}^{\ast }(X)$.

\bigskip

\textbf{Proposition 2.2 }(see [23])\textbf{.} \textit{If }$H$\textit{\ and }$%
K$\textit{\ are two nonempty subsets of the metric space }$(X,d)$\textit{,
then }%
\begin{equation*}
h(H,K)=h(\overline{H},\overline{K})\text{.}
\end{equation*}

\bigskip

\textbf{Proposition 2.3 }(see [23])\textbf{.} \textit{If }$(H_{i})_{i\in I}$%
\textit{\ and }$(K_{i})_{i\in I}$\textit{\ are two families of nonempty
subsets of the metric space }$(X,d)$\textit{, then}

\begin{equation*}
h(\underset{i\in I}{\cup }H_{i},\underset{i\in I}{\cup }K_{i})=h(\overline{%
\underset{i\in I}{\cup }H_{i}},\overline{\underset{i\in I}{\cup }K_{i}})\leq 
\underset{i\in I}{\sup }h(H_{i},K_{i})\text{.}
\end{equation*}

\bigskip

\textbf{Theorem 2.4 }(see [23])\textbf{.}\textit{\ If the metric space }$%
(X,d)$ \textit{is complete,} \textit{then} $(\mathcal{K}(X),h)$\textit{\ is
a complete metric space.}

\bigskip

\textbf{Definition 2.5.} \textit{For a metric space }$(X,d)$\textit{, we
consider on} $\mathcal{P}^{\ast }(X)$ \textit{the function }$\delta :%
\mathcal{P}^{\ast }(X)\times \mathcal{P}^{\ast }(X)\rightarrow \lbrack
0,+\infty ]$ \textit{defined by}

\begin{equation*}
\delta (A,B)=\underset{x\in A,y\in B}{\sup }d(x,y)\text{,}
\end{equation*}%
\textit{for all} $A,B\in \mathcal{P}^{\ast }(X)$.

\bigskip

\textbf{Remark 2.6.} \textit{For every }$A,B\in \mathcal{P}^{\ast }(X)$%
\textit{\ we have}%
\begin{equation*}
h(A,B)\leq \delta (A,B)\text{.}
\end{equation*}

\bigskip

\textbf{Proposition 2.7.} \textit{Let }$(X,d)$\textit{\ be a complete metric
space, }$(Y_{n})_{n\in \mathbb{N}}\subseteq \mathcal{K}(X)$ \textit{and }$Y$ 
\textit{a closed subset of }$X$\textit{\ such that }$\underset{n\rightarrow
\infty }{\lim }h(Y_{n},Y)=0$\textit{. Then }$Y\in \mathcal{K}(X)$\textit{.}

\textit{Proof. }It is enough to prove that $Y$ is precompact. To this aim,
let us note that for each $\varepsilon >0$ there exists $n_{\varepsilon }\in 
\mathbb{N}$ such that $h(Y_{n_{\varepsilon }},Y)<\frac{\varepsilon }{2}$, so 
$Y\subseteq N_{\frac{\varepsilon }{2}}(Y_{n_{\varepsilon }})$. Since $%
Y_{n_{\varepsilon }}\in \mathcal{K}(X)$\textit{\ }there exist $%
x_{1},...,x_{m}\in X$ such that $Y_{n_{\varepsilon }}\subseteq \underset{i=1}%
{\overset{m}{\cup }}B(x_{i},\frac{\varepsilon }{2})$ and therefore $%
Y\subseteq \underset{i=1}{\overset{m}{\cup }}B(x_{i},\varepsilon )$. $%
\square $

\bigskip

\textbf{Proposition 2.8.} \textit{Let }$(X,d)$\textit{\ be a complete metric
space, }$(Y_{n})_{n\in \mathbb{N}}\subseteq \mathcal{K}(X)$ \textit{and }$%
Y\in \mathcal{K}(X)$\textit{\ such that }$\underset{n\rightarrow \infty }{%
\lim }h(Y_{n},Y)=0$\textit{. Then }$H\overset{def}{=}Y\cup (\overset{\infty }%
{\underset{n=0}{\cup }}Y_{n})\in \mathcal{K}(X)$\textit{.}

\textit{Proof. }First of all we prove that $H$ is a closed subset of $X$.

Indeed, if $x\in \overline{H}$, then there exists a sequence $(x_{k})_{k\in 
\mathbb{N}}\subseteq H$ such that $\underset{k\rightarrow \infty }{\lim }%
x_{k}=x$.

If $\{k\in \mathbb{N\mid }x_{k}\in Y\}$ is infinite, then there exists a
subsequence $(x_{k_{p}})_{p\in \mathbb{N}}$ of $(x_{k})_{k\in \mathbb{N}}$
such that $x_{k_{p}}\in Y$ for every $p\in \mathbb{N}$. Since $Y\in \mathcal{%
K}(X)$ there exists a subsequence $(x_{k_{p_{q}}})_{q\in \mathbb{N}}$ of $%
(x_{k_{p}})_{p\in \mathbb{N}}$ and $y\in Y$ such that $\underset{%
q\rightarrow \infty }{\lim }x_{k_{p_{q}}}=y$. Consequently, as $\underset{%
q\rightarrow \infty }{\lim }x_{k_{p_{q}}}=x$, we conclude that $x=y\in
Y\subseteq H$.

If there exists $n_{0}\in \mathbb{N}$ such that $\{k\in \mathbb{N\mid }%
x_{k}\in Y_{n_{0}}\}$ is infinite, a similar argument shows that $x\in H$.

If none of the above described two cases is valid, then there exist an
increasing sequence $(k_{p})_{p\in \mathbb{N}}\subseteq \mathbb{N}$, $%
x_{k_{p}}\in Y_{k_{p}}$ and $y_{k_{p}}\in Y$ such that 
\begin{equation*}
d(x_{k_{p}},y_{k_{p}})<h(Y_{k_{p}},Y)+\frac{1}{p}\text{.}
\end{equation*}%
Since $Y\in \mathcal{K}(X)$\textit{\ }there exists $(y_{k_{p_{q}}})_{q\in 
\mathbb{N}}$ a subsequence of $(y_{k_{p}})_{p\in \mathbb{N}}$ and $y\in Y$
such that $\underset{q\rightarrow \infty }{\lim }y_{k_{p_{q}}}=y$. As 
\begin{equation*}
d(x_{k_{p_{q}}},y)<d(x_{k_{p_{q}}},y_{k_{p_{q}}})+d(y_{k_{p_{q}}},y)\leq
h(Y_{k_{p_{q}}},Y)+\frac{1}{p_{q}}+d(y_{k_{p_{q}}},y)
\end{equation*}%
and 
\begin{equation*}
\underset{q\rightarrow \infty }{\lim }h(Y_{k_{p_{q}}},Y)=\underset{%
q\rightarrow \infty }{\lim }\frac{1}{p_{q}}=\underset{q\rightarrow \infty }{%
\lim }d(y_{k_{p_{q}}},y)=0\text{,}
\end{equation*}%
we infer that $\underset{q\rightarrow \infty }{\lim }x_{k_{p_{q}}}=y$.
Consequently, as $\underset{q\rightarrow \infty }{\lim }x_{k_{p_{q}}}=x$, we
get $x=y\in Y\subseteq H$.

Now we prove that 
\begin{equation*}
\underset{m\rightarrow \infty }{\lim }h(\underset{i=0}{\overset{m}{\cup }}%
Y_{i},H)=0\text{.}
\end{equation*}

Indeed%
\begin{equation*}
h(\underset{i=0}{\overset{m}{\cup }}Y_{i},H)=h((\underset{i=0}{\overset{m}{%
\cup }}Y_{i})\cup (\underset{i=m+1}{\overset{\infty }{\cup }}Y_{m})\cup
Y_{m},(\underset{i=0}{\overset{m}{\cup }}Y_{i})\cup (\underset{i=m+1}{%
\overset{\infty }{\cup }}Y_{i})\cup Y))\overset{\text{Proposition 2.3}}{\leq 
}
\end{equation*}%
\begin{equation*}
\leq \sup \{h(Y_{m},Y),h(Y_{m},Y_{m+1}),h(Y_{m},Y_{m+2}),...\}
\end{equation*}%
for every $m\in \mathbb{N}$. As $\underset{m\rightarrow \infty }{\lim }%
h(Y_{m},Y)=0$, we conclude that $\underset{m\rightarrow \infty }{\lim }h(%
\underset{i=0}{\overset{m}{\cup }}Y_{i},H)=0$.

Because $\underset{i=0}{\overset{m}{\cup }}Y_{i}\in \mathcal{K}(X)$ for
every $m\in \mathbb{N}$ and $H$ is closed, using Proposition 2.7, we obtain
that $H$ is compact. $\square $

\bigskip

\textbf{3. The main results}

\bigskip

\textbf{Definition 3.1.\ }\textit{An iterated function system consisting of
convex contractions on a complete metric space }$(X,d)$\textit{\ is given by
a finite family of continuous functions }$(f_{i})_{i\in I}$\textit{, }$%
f_{i}:X\rightarrow X$,\textit{\ such that} \textit{for every }$i,j\in I$%
\textit{\ there exist }$a_{ij},b_{ij},c_{ij}\in \lbrack 0,\infty )$\textit{\
satisfying the following two properties:}

$\alpha $) $a_{ij}+b_{ij}+c_{ij}\overset{def}{=}d_{ij}$ \textit{and} $%
\underset{i,j\in I}{\max }d_{ij}\overset{def}{=}d<1$;

$\beta $) 
\begin{equation*}
d((f_{i}\circ f_{j})(x),(f_{i}\circ f_{j})(y))\leq
a_{ij}d(x,y)+b_{ij}d(f_{i}(x),f_{i}(y))+c_{ij}d(f_{j}(x),f_{j}(y))
\end{equation*}%
\textit{for every }$i,j\in I$\textit{\ and every }$x,y\in X$\textit{.}

\textit{We denote such a system by }%
\begin{equation*}
\mathcal{S}=((X,d),(f_{i})_{i\in I})\text{\textit{.}}
\end{equation*}

\bigskip

One can associate to the system $\mathcal{S}$\ the function $F_{\mathcal{S}}:%
\mathcal{K}(X)\rightarrow \mathcal{K}(X)$\ given by%
\begin{equation*}
F_{\mathcal{S}}(B)=\underset{i\in I}{\cup }f_{i}(B)
\end{equation*}%
for all $B\in \mathcal{K}(X)$.

\bigskip

\textbf{Theorem 3.2.} \textit{Let} $\mathcal{S}=((X,d),(f_{i})_{i\in I})$ 
\textit{be an iterated function system consisting of convex contractions.
Then:}

i) \textit{There exists a unique }$A\in \mathcal{K}(X)$\textit{\ such that }%
\begin{equation*}
\underset{n\rightarrow \infty }{\lim }h(F_{\mathcal{S}}^{[n]}(B),A)=0\text{,}
\end{equation*}%
\textit{for every }$B\in \mathcal{K}(X)$\textit{.}

ii) \textit{For each }$\omega \in \Lambda (I)$\textit{\ there exists }$%
a_{\omega }\in X$\textit{\ such that} 
\begin{equation*}
\underset{n\rightarrow \infty }{\lim }h(f_{[\omega ]_{n}}(B),\{a_{\omega
}\})=0\text{,}
\end{equation*}%
\textit{for every }$B\in \mathcal{K}(X)$\textit{.}

\textit{Moreover}%
\begin{equation*}
\underset{n\rightarrow \infty }{\lim }\underset{\omega \in \Lambda (I)}{\sup 
}h(f_{[\omega ]_{n}}(B),\{a_{\omega }\})=0
\end{equation*}%
\textit{for every }$B\in \mathcal{K}(X)$\textit{.}

iii)%
\begin{equation*}
A=\overline{\{a_{\omega }\mid \omega \in \Lambda (I)\}}.
\end{equation*}

iv) \textit{For every} $(Y_{n})_{n\in \mathbb{N}}\subseteq \mathcal{K}(X)$%
\textit{\ and\ }$Y\in \mathcal{K}(X)$, \textit{the following implication is
valid: }%
\begin{equation*}
\underset{n\rightarrow \infty }{\lim }h(Y_{n},Y)=0\Rightarrow \underset{%
n\rightarrow \infty }{\lim }h(F_{\mathcal{S}}(Y_{n}),F_{\mathcal{S}}(Y))=0%
\text{\textit{.}}
\end{equation*}

v)\textit{\ }$A$\textit{\ is the unique fixed point of }$F_{\mathcal{S}}$.

\textit{Proof}.

i) For fixed $Y,Z\in \mathcal{K}(X)$ we define%
\begin{equation*}
x_{n}(Y,Z)=\underset{\omega \in \Lambda _{n}(I)}{\sup }\delta (f_{\omega
}(Y),f_{\omega }(Z))
\end{equation*}%
and%
\begin{equation*}
y_{n}(Y,Z)=\max \{x_{n-1}(Y,Z),x_{n}(Y,Z)\}
\end{equation*}%
for every $n\in \mathbb{N}^{\ast }$. For the sake of simplicity we will
denote $x_{n}(Y,Z)$ by $x_{n}$ and $y_{n}(Y,Z)$ by $y_{n}$.

We claim that the sequence $(y_{n})_{n\in \mathbb{N}^{\ast }}$ is decreasing.

Indeed, for $n\in \mathbb{N}^{\ast }$ and $\omega \in \Lambda _{n+1}(I)$
there exist $i,j\in I$ and $\omega _{0}\in \Lambda _{n-1}(I)$ such that $%
\omega =ij\omega _{0}$. Then, for $y\in Y$ and $z\in Z$, we have%
\begin{equation*}
d(f_{\omega }(y),f_{\omega }(z))=d(f_{ij\omega _{0}}(y),f_{ij\omega
_{0}}(z))\leq
\end{equation*}%
\begin{equation*}
\leq a_{ij}d(f_{\omega _{0}}(y),f_{\omega _{0}}(z))+b_{ij}d(f_{i\omega
_{0}}(y),f_{i\omega _{0}}(z))+c_{ij}d(f_{j\omega _{0}}(y),f_{j\omega
_{0}}(z))\leq
\end{equation*}%
\begin{equation*}
\leq a_{ij}x_{n-1}+b_{ij}x_{n}+c_{ij}x_{n}\leq
a_{ij}x_{n-1}+(b_{ij}+c_{ij})x_{n}\leq
\end{equation*}%
\begin{equation*}
\leq d_{ij}\max \{x_{n-1},x_{n}\}=d_{ij}y_{n}\leq dy_{n}<y_{n}\text{,}
\end{equation*}%
so%
\begin{equation}
x_{n+1}=\underset{\omega \in \Lambda _{n+1}(I)}{\sup }\delta (f_{\omega
}(Y),f_{\omega }(Z))\leq dy_{n}<y_{n}\text{.}  \tag{1}
\end{equation}%
As 
\begin{equation}
x_{n}\leq \max \{x_{n-1},x_{n}\}=y_{n}\text{,}  \tag{2}
\end{equation}%
we get%
\begin{equation*}
y_{n+1}=\max \{x_{n},x_{n+1}\}\overset{(1)\text{ and }(2)}{\leq }y_{n}\text{.%
}
\end{equation*}

Therefore we have%
\begin{equation*}
y_{n+2}=\max \{x_{n+1},x_{n+2}\}\overset{(1)}{\leq }\max
\{dy_{n},dy_{n+1}\}=dy_{n}
\end{equation*}%
and consequently%
\begin{equation*}
y_{2n-1}\leq d^{n-1}y_{1}
\end{equation*}%
and%
\begin{equation*}
y_{2n}\leq d^{n-1}y_{1}
\end{equation*}%
for every $n\in \mathbb{N}^{\ast }$.

Thus the series $\underset{n\in \mathbb{N}^{\ast }}{\sum }y_{n}$ is
convergent, so the series $\underset{n\in \mathbb{N}^{\ast }}{\sum }x_{n}$
is convergent (see $(2)$ and use the comparison test) and consequently $%
\underset{n\rightarrow \infty }{\lim }x_{n}=0$. Hence, as%
\begin{equation*}
h(F_{\mathcal{S}}^{[n]}(Y),F_{\mathcal{S}}^{[n]}(Z))=h(\underset{\omega \in
\Lambda _{n}(I)}{\cup }f_{\omega }(Y),\underset{\omega \in \Lambda _{n}(I)}{%
\cup }f_{\omega }(Z))\overset{\text{Proposition 2.3}}{\leq }
\end{equation*}%
\begin{equation}
\leq \underset{\omega \in \Lambda _{n}(I)}{\sup }h(f_{\omega }(Y),f_{\omega
}(Z))\overset{\text{Remark 2.6}}{\leq }\underset{\omega \in \Lambda _{n}(I)}{%
\sup }\delta (f_{\omega }(Y),f_{\omega }(Z))=x_{n}  \tag{3}
\end{equation}%
for every $n\in \mathbb{N}^{\ast }$, we get that%
\begin{equation}
\underset{n\rightarrow \infty }{\lim }h(F_{\mathcal{S}}^{[n]}(Y),F_{\mathcal{%
S}}^{[n]}(Z))=0\text{.}  \tag{4}
\end{equation}

In particular, for each $Y\in \mathcal{K}(X)$, considering $Z=F_{\mathcal{S}%
}(Y)\in \mathcal{K}(X)$ and taking into account the comparison test and $(3)$%
, we infer that the series $\underset{n\in \mathbb{N}^{\ast }}{\sum }h(F_{%
\mathcal{S}}^{[n+1]}(Y),F_{\mathcal{S}}^{[n]}(Y))$ is convergent. Thus we
conclude that the sequence $(F_{\mathcal{S}}^{[n+1]}(Y))_{n\in \mathbb{N}%
^{\ast }}$ is Cauchy and, as $(\mathcal{K}(X),h)$ is complete (see Theorem
2.4), there exists $A_{Y}\in \mathcal{K}(X)$ such that%
\begin{equation}
\underset{n\rightarrow \infty }{\lim }h(F_{\mathcal{S}}^{[n]}(Y),A_{Y})=0%
\text{.}  \tag{5}
\end{equation}

In the same manner we can prove that if $Z\in \mathcal{K}(X)$, then%
\begin{equation}
\underset{n\rightarrow \infty }{\lim }h(F_{\mathcal{S}}^{[n]}(Z),A_{Z})=0%
\text{.}  \tag{6}
\end{equation}

From $(4)$, $(5)$ and $(6)$ we obtain that $A_{Y}=A_{Z}\overset{def}{=}A$
for every $Y,Z\in \mathcal{K}(X)$. Thus 
\begin{equation*}
\underset{n\rightarrow \infty }{\lim }h(F_{\mathcal{S}}^{[n]}(B),A)=0\text{,}
\end{equation*}%
for every\textit{\ }$B\in \mathcal{K}(X)$\textit{.}

ii) For $\omega \in \Lambda (I)$ and $Y,Z\in \mathcal{K}(X)$ we have%
\begin{equation*}
h(f_{[\omega ]_{n}}(Y),f_{[\omega ]_{n}}(Z))\overset{\text{Remark 2.6}}{\leq 
}\delta (f_{[\omega ]_{n}}(Y),f_{[\omega ]_{n}}(Z))\leq \underset{\omega \in
\Lambda _{n}(I)}{\sup }\delta (f_{\omega }(Y),f_{\omega }(Z))=x_{n}
\end{equation*}%
for every $n\in \mathbb{N}^{\ast }$, so, as $\underset{n\rightarrow \infty }{%
\lim }x_{n}=0$, we deduce that%
\begin{equation}
\underset{n\rightarrow \infty }{\lim }\delta (f_{[\omega
]_{n}}(Y),f_{[\omega ]_{n}}(Z))=\underset{n\rightarrow \infty }{\lim }%
h(f_{[\omega ]_{n}}(Y),f_{[\omega ]_{n}}(Z))=0\text{.}  \tag{7}
\end{equation}

For $Y\in \mathcal{K}(X)$ we have%
\begin{equation*}
h(f_{[\omega ]_{n}}(Y),f_{[\omega ]_{n+1}}(Y))\overset{\text{Remark 2.6}}{%
\leq }\delta (f_{[\omega ]_{n}}(Y),f_{[\omega ]_{n+1}}(Y))\leq
\end{equation*}%
\begin{equation*}
\leq \delta (f_{[\omega ]_{n}}(Y),f_{[\omega ]_{n}}(F_{\mathcal{S}}(Y)))\leq
x_{n}(Y,F_{\mathcal{S}}(Y))
\end{equation*}%
for each $n\in \mathbb{N}^{\ast }$, hence, since -as we have seen in the
proof of 1)- the series $\underset{n\in \mathbb{N}^{\ast }}{\sum }x_{n}(Y,F_{%
\mathcal{S}}(Y))$ is convergent, using the comparison criterion, we infer
that the series $\underset{n\in \mathbb{N}^{\ast }}{\sum }h(f_{[\omega
]_{n}}(Y),f_{[\omega ]_{n+1}}(Y))$ is convergent. Thus we conclude that the
sequence $(f_{[\omega ]_{n}}(Y))_{n\in \mathbb{N}^{\ast }}$ is Cauchy and
as, $(\mathcal{K}(X),h)$ is complete (see Theorem 2.4), there exists $%
A_{\omega }(Y)\in \mathcal{K}(X)$ such that%
\begin{equation}
\underset{n\rightarrow \infty }{\lim }h(f_{[\omega ]_{n}}(Y),A_{\omega
}(Y))=0\text{.}  \tag{8}
\end{equation}

In the same manner we can prove that if $Z\in \mathcal{K}(X)$, then there
exists $A_{\omega }(Z)\in \mathcal{K}(X)$ such that%
\begin{equation}
\underset{n\rightarrow \infty }{\lim }h(f_{[\omega ]_{n}}(Z),A_{\omega
}(Z))=0\text{.}  \tag{9}
\end{equation}

From $(7)$, $(8)$ and $(9)$ we obtain that $A_{\omega }(Y)=A_{\omega }(Z)%
\overset{def}{=}A_{\omega }$ for each $Y,Z\in \mathcal{K}(X)$. Thus 
\begin{equation}
\underset{n\rightarrow \infty }{\lim }h(f_{[\omega ]_{n}}(B),A_{\omega })=0%
\text{,}  \tag{10}
\end{equation}%
for each\textit{\ }$B\in \mathcal{K}(X)$\textit{.}

Since 
\begin{equation}
\underset{n\rightarrow \infty }{\lim }diam(f_{[\omega ]_{n}}(B))=0  \tag{11}
\end{equation}%
for each\textit{\ }$B\in \mathcal{K}(X)$ (see $(7)$ for $Y=Z=B$), we get
that 
\begin{equation}
diam(A_{\omega })=0\text{.}  \tag{12}
\end{equation}

Indeed, using $(10)$ and $(11)$, we infer that for each $\varepsilon >0$
there exists $n_{\varepsilon }\in \mathbb{N}^{\ast }$ such that 
\begin{equation*}
diam(f_{[\omega ]_{n_{\varepsilon }}}(B))<\varepsilon \text{ and }%
h(f_{[\omega ]_{n_{\varepsilon }}}(B),A_{\omega })<\varepsilon .
\end{equation*}
Therefore there exists $\eta _{0}\in (0,\varepsilon )$ such that 
\begin{equation*}
A_{\omega }\subseteq N_{\eta _{0}}(f_{[\omega ]_{n_{\varepsilon }}}(B))\text{%
,}
\end{equation*}
so 
\begin{equation*}
diam(A_{\omega })\leq 2\eta _{0}+diam(f_{[\omega ]_{n_{\varepsilon
}}}(B))<3\varepsilon \text{.}
\end{equation*}%
As $\varepsilon $ was arbitrarily chosen, we conclude that $diam(A_{\omega
})=0$.

From $(12)$ we conclude that there exists $a_{\omega }\in X$\textit{\ }such
that $A_{\omega }=\{a_{\omega }\}$ and, from $(10)$, we get 
\begin{equation*}
\underset{n\rightarrow \infty }{\lim }h(f_{[\omega ]_{n}}(B),\{a_{\omega
}\})=0\text{,}
\end{equation*}%
for each\textit{\ }$B\in \mathcal{K}(X)$\textit{.}

Note that the above limit is uniform with respect to $\omega \in \Lambda (I)$%
, i.e.%
\begin{equation*}
\underset{n\rightarrow \infty }{\lim }\underset{\omega \in \Lambda (I)}{\sup 
}h(f_{[\omega ]_{n}}(B),\{a_{\omega }\})=0\text{.}
\end{equation*}

Indeed, 
\begin{equation*}
h(f_{[\omega ]_{n}}(B),\{a_{\omega }\})\leq \underset{k=n}{\overset{m}{\sum }%
}h(f_{[\omega ]_{k}}(B),f_{[\omega ]_{k+1}}(B))+h(f_{[\omega
]_{m+1}}(B),\{a_{\omega }\})
\end{equation*}%
for every $m,n\in \mathbb{N}$, $m\geq n$. By passing to limit as $%
m\rightarrow \infty $, we get 
\begin{equation*}
h(f_{[\omega ]_{n}}(B),\{a_{\omega }\})\leq \underset{k\geq n}{\sum }%
h(f_{[\omega ]_{k}}(B),f_{[\omega ]_{k+1}}(B))\overset{\text{Remark 2.6}}{%
\leq }
\end{equation*}%
\begin{equation*}
\leq \underset{k\geq n}{\sum }\delta (f_{[\omega ]_{k}}(B),f_{[\omega
]_{k}}(F_{\mathcal{S}}(B)))=\underset{k\geq n}{\sum }x_{k}(B,F_{\mathcal{S}%
}(B))
\end{equation*}%
for every $\omega \in \Lambda (I)$ and every $n\in \mathbb{N}$, so 
\begin{equation*}
\underset{\omega \in \Lambda (I)}{\sup }h(f_{[\omega ]_{n}}(B),\{a_{\omega
}\})\leq \underset{k\geq n}{\sum }x_{k}(B,F_{\mathcal{S}}(B))
\end{equation*}%
for every $n\in \mathbb{N}$. As the series $\underset{n}{\sum }x_{n}(B,F_{%
\mathcal{S}}(B))$ is convergent, we conclude that $\underset{n\rightarrow
\infty }{\lim }\underset{\omega \in \Lambda (I)}{\sup }h(f_{[\omega
]_{n}}(B),\{a_{\omega }\})=0$.

iii) As%
\begin{equation*}
h(F_{\mathcal{S}}^{[n]}(B),\{a_{\omega }\mid \omega \in \Lambda (I)\})=
\end{equation*}%
\begin{equation*}
=h(\underset{\omega \in \Lambda _{n}\left( I\right) }{\cup }\underset{\alpha
\in \Lambda \left( I\right) }{\cup }f_{[\omega \alpha ]_{n}}(B),\underset{%
\omega \in \Lambda _{n}\left( I\right) }{\cup }\underset{\alpha \in \Lambda
\left( I\right) }{\cup }\{a_{\omega \alpha }\mid \alpha \in \Lambda (I)\})%
\overset{\text{Proposition 2.3}}{\leq }
\end{equation*}%
\begin{equation*}
\leq \underset{\omega \in \Lambda _{n}\left( I\right) }{\sup }\underset{%
\alpha \in \Lambda \left( I\right) }{\sup }h(f_{\omega }(B),\{a_{\omega
\alpha }\})\text{,}
\end{equation*}%
we have%
\begin{equation*}
h(A,\{a_{\omega }\mid \omega \in \Lambda (I)\})\leq h(A,F_{\mathcal{S}%
}^{[n]}(B))+h(F_{\mathcal{S}}^{[n]}(B),\{a_{\omega }\mid \omega \in \Lambda
(I)\})\leq
\end{equation*}%
\begin{equation}
\leq h(A,F_{\mathcal{S}}^{[n]}(B))+\underset{\omega \in \Lambda _{n}\left(
I\right) }{\sup }\underset{\alpha \in \Lambda \left( I\right) }{\sup }%
h(f_{\omega }(B),\{a_{\omega \alpha }\})  \tag{13}
\end{equation}%
for all $n\in \mathbb{N}^{\ast }$ and $B\in \mathcal{K}(X)$\textit{.}

Since 
\begin{equation*}
\underset{n\rightarrow \infty }{\lim }h(F_{\mathcal{S}}^{[n]}(B),A)=0
\end{equation*}%
(see i)) and 
\begin{equation*}
\underset{n\rightarrow \infty }{\lim }\underset{\omega \in \Lambda
_{n}\left( I\right) }{\sup }\underset{\alpha \in \Lambda \left( I\right) }{%
\sup }h(f_{\omega }(B),\{a_{\omega \alpha }\})=0
\end{equation*}%
(see ii)), by passing to limit in $(13)$, we obtain that 
\begin{equation*}
h(A,\{a_{\omega }\mid \omega \in \Lambda (I)\})=0\text{,}
\end{equation*}%
i.e. 
\begin{equation*}
h(A,\overline{\{a_{\omega }\mid \omega \in \Lambda (I)\}})=0
\end{equation*}%
(see Proposition 2.2).

Hence 
\begin{equation*}
A=\overline{\{a_{\omega }\mid \omega \in \Lambda (I)\}}\text{.}
\end{equation*}

iv) Let us consider\textit{\ }$(Y_{n})_{n\in \mathbb{N}}\subseteq \mathcal{K}%
(X)$\textit{\ }and\textit{\ }$Y\in \mathcal{K}(X)$ such that\textit{\ }%
\begin{equation*}
\underset{n\rightarrow \infty }{\lim }h(Y_{n},Y)=0\text{.}
\end{equation*}

Using Proposition 2.8 we conclude that 
\begin{equation*}
H\overset{def}{=}Y\cup (\overset{\infty }{\underset{n=0}{\cup }}Y_{n})\in 
\mathcal{K}(X)\text{.}
\end{equation*}
Hence, as the functions $f_{i}$ are continuous, they are uniformly
continuous on $H$, so for each $\varepsilon >0$ there exists $\delta
_{\varepsilon }>0$ such that 
\begin{equation*}
d(f_{i}(x),f_{i}(y))<\frac{\varepsilon }{2}
\end{equation*}
for every $i\in I$ and every $x,y\in H$ such that $d(x,y)<\delta
_{\varepsilon }$.

For each $\varepsilon >0$ there exists $n_{\varepsilon }\in \mathbb{N}$ such
that 
\begin{equation*}
h(Y_{n},Y)<\frac{\delta _{\varepsilon }}{2}
\end{equation*}
for every $n\in \mathbb{N}$, $n\geq n_{\varepsilon }$.

Let us consider $i\in I$ and $n\in \mathbb{N}$, $n\geq n_{\varepsilon }$.

Since for every $x\in Y_{n}\subseteq H$ there exists $y\in Y\subseteq H$
such that 
\begin{equation*}
d(x,y)<d(x,Y)+\frac{\delta _{\varepsilon }}{2}\text{,}
\end{equation*}%
we get that 
\begin{equation*}
d(x,y)<h(Y_{n},Y)+\frac{\delta _{\varepsilon }}{2}<\frac{\delta
_{\varepsilon }}{2}+\frac{\delta _{\varepsilon }}{2}=\delta _{\varepsilon }%
\text{,}
\end{equation*}%
so 
\begin{equation*}
d(f_{i}(x),f_{i}(Y))\leq d(f_{i}(x),f_{i}(y))<\frac{\varepsilon }{2}\text{.}
\end{equation*}%
Consequently 
\begin{equation*}
d(f_{i}(Y_{n}),f_{i}(Y))\leq \frac{\varepsilon }{2}\text{.}
\end{equation*}%
In the same manner, one can prove that 
\begin{equation*}
d(f_{i}(Y),f_{i}(Y_{n}))\leq \frac{\varepsilon }{2}\text{,}
\end{equation*}%
so 
\begin{equation*}
h(f_{i}(Y),f_{i}(Y_{n}))\leq \frac{\varepsilon }{2}\text{.}
\end{equation*}%
Hence 
\begin{equation*}
h(F_{\mathcal{S}}(Y_{n}),F_{\mathcal{S}}(Y))=h(\underset{i\in I}{\cup }%
f_{i}(Y_{n}),\underset{i\in I}{\cup }f_{i}(Y))\overset{\text{Proposition 2.3}%
}{\leq }
\end{equation*}%
\begin{equation*}
\leq \underset{i\in I}{\max }h(f_{i}(Y_{n}),f_{i}(Y))\leq \frac{\varepsilon 
}{2}<\varepsilon \text{.}
\end{equation*}

Thus for each $\varepsilon >0$ there exists $n_{\varepsilon }\in \mathbb{N}$
such that $h(F_{\mathcal{S}}(Y_{n}),F_{\mathcal{S}}(Y))<\varepsilon $ for
every $n\in \mathbb{N}$, $n\geq n_{\varepsilon }$, i.e. 
\begin{equation*}
\underset{n\rightarrow \infty }{\lim }h(F_{\mathcal{S}}(Y_{n}),F_{\mathcal{S}%
}(Y))=0\text{\textit{.}}
\end{equation*}

v) Since $\underset{n\rightarrow \infty }{\lim }h(F_{\mathcal{S}%
}^{[n]}(A),A)=0$ (see i) for $B=A$), using 4) for $Y_{n}=F_{\mathcal{S}%
}^{[n]}(A)\in \mathcal{K}(X)$ and $Y=A\in \mathcal{K}(X)$, we obtain that 
\begin{equation}
\underset{n\rightarrow \infty }{\lim }h(F_{\mathcal{S}}^{[n+1]}(A),F_{%
\mathcal{S}}(A))=0\text{.}  \tag{14}
\end{equation}

Using i), for $B=F_{\mathcal{S}}(A)$, we infer that%
\begin{equation}
\underset{n\rightarrow \infty }{\lim }h(F_{\mathcal{S}}^{[n+1]}(A),A)=0\text{%
.}  \tag{15}
\end{equation}

From $(14)$ and $(15)$ we conclude that 
\begin{equation*}
F_{\mathcal{S}}(A)=A\text{.}
\end{equation*}

Moreover, if for some $A_{1}\in \mathcal{K}(X)$ we have $F_{\mathcal{S}%
}(A_{1})=A_{1}$, then $F_{\mathcal{S}}^{[n]}(A_{1})=A_{1}$ for each $n\in 
\mathbb{N}$, so $\underset{n\rightarrow \infty }{\lim }h(F_{\mathcal{S}%
}^{[n]}(A_{1}),A_{1})=0$. Since, according to i), we have $\underset{%
n\rightarrow \infty }{\lim }h(F_{\mathcal{S}}^{[n]}(A_{1}),A)=0$, we
conclude that $A=A_{1}$. $\square $

\bigskip

Let us note that, concerning the speed of convergence of the sequence $(F_{%
\mathcal{S}}^{[n]}(B))_{n\in \mathbb{N}}$, where $B\in \mathcal{K}(X)$, we
have (from the proof of i), the following inequality:

\begin{equation*}
h(F_{\mathcal{S}}^{[n]}(B),A)\leq \frac{d^{[\frac{n}{2}]}}{1-d}(x_{0}(B,F_{%
\mathcal{S}}(B))+x_{1}(B,F_{\mathcal{S}}(B)))\text{,}
\end{equation*}%
for every $n\in \mathbb{N}$.

\bigskip

\textbf{Remark 3.3.} \textit{By taking in the above Theorem\ a set }$I$%
\textit{\ with one element, we get that }$A$ \textit{has exactly one element
which is the fixed point of the convex contraction that can be approximated
by means of Picard iteration. Consequently we obtain Istr\u{a}\c{t}escu's
fixed point theorem for convex contractions.}

\bigskip

\textbf{Proposition 3.4}. \textit{Let} $\mathcal{S}=((X,d),(f_{i})_{i\in I})$
\textit{be an iterated function system consisting of convex contractions.
Then, in the framework of Theorem 3.2, we have }%
\begin{equation*}
\underset{n\rightarrow \infty }{\lim }diam(A_{[\omega ]_{n}})=0
\end{equation*}%
\textit{for every} $\omega \in \Lambda (I)$.

\textit{Proof}. Take $B=A$ in $(11)$ from the proof of Theorem 3.2. $\square 
$

\bigskip

\textbf{Proposition 3.5.} \textit{Let} $\mathcal{S}=((X,d),(f_{i})_{i\in I})$
\textit{be an iterated function system consisting of convex contractions.
Then, in the framework of Theorem 3.2, we have }%
\begin{equation*}
\underset{n\in \mathbb{N}}{\cap }A_{[\omega ]_{n}}=\{a_{\omega }\}
\end{equation*}%
\textit{for every} $\omega \in \Lambda (I)$.

\textit{Proof}. From $F_{\mathcal{S}}(A)=A$ we infer that 
\begin{equation*}
A_{[\omega ]_{n+1}}\subseteq A_{[\omega ]_{n}}
\end{equation*}%
for every $n\in \mathbb{N}$. Then 
\begin{equation*}
\underset{n\rightarrow \infty }{\lim }h(A_{[\omega ]_{n}},\underset{n\in 
\mathbb{N}}{\cap }A_{[\omega ]_{n}})=0
\end{equation*}%
(see Theorem 1.14 from [23]) and taking into account Theorem 3.2, ii), we
conclude that $\underset{n\in \mathbb{N}}{\cap }A_{[\omega
]_{n}}=\{a_{\omega }\}$. $\square $

\bigskip

Using the above two Propositions, the same arguments as the ones used in the
proof of Theorem 4.1\ from [18] give us the following:

\bigskip

\textbf{Theorem 3.6.} \textit{Let} $\mathcal{S}=((X,d),(f_{i})_{i\in I})$ 
\textit{be an iterated function system consisting of convex contractions.
Then, in the framework of Theorem 3.2, the function }$\pi :\Lambda
(I)\rightarrow A$\textit{\ defined by }%
\begin{equation*}
\pi (\omega )=a_{\omega }\text{,}
\end{equation*}%
\textit{for every\ }$\omega \in \Lambda $\textit{, which is called the
canonical projection from }$\Lambda (I)$ \textit{to} $A$,\textit{\ has the
following properties:}

1)\textit{\ it is continuous;}

2) \textit{it\ is onto;}

3) 
\begin{equation*}
\pi \circ F_{i}=f_{i}\circ \pi \text{,}
\end{equation*}%
\textit{\ for every }$i\in I$.

\bigskip

\textbf{References}

\bigskip

[1] J. Andres and M. Rypka, Multivalued fractals and hyperfractals,
Internat. J. Bifur. Chaos Appl. Sci. Engrg., \textbf{22} (2012), DOI
10.1142/S02181127412500095.

[2] M. Boriceanu, M. Bota and A. Petru\c{s}el, Multivalued fractals in $b$%
-metric spaces, Cent. Eur. J. Math., \textbf{8} (2010), 367-377.

[3] C. Chifu and A. Petru\c{s}el, Multivalued fractals and generalized
multivalued contractions, Chaos Solitons Fractals, \textbf{36} (2008),
203-210.

[4] D. Dumitru, Generalized iterated function systems containing Meir-Keeler
functions, An. Univ. Bucur., Mat., \textbf{58} (2009), 109-121.

[5] V. Istr\u{a}\c{t}escu, Some fixed point theorems for convex contraction
mappings and convex nonexpansive mappings (I), Libertas Math., \textbf{1}
(1981), 151-164.

[6] V. Istr\u{a}\c{t}escu, Some fixed point theorems for convex contraction
mappings and mappings with convex diminishing diameters - I, Annali di Mat.
Pura Appl., \textbf{130} (1982), 89--104.

[7] V. Istr\u{a}\c{t}escu,\ Some fixed point theorems for convex contraction
mappings and mappings with convex diminishing diameters, II, Annali di Mat.
Pura Appl., \textbf{134} (1983), 327-362.

[8] V. Ghorbanian, S. Rezapour and N. Shahzad, Some ordered fixed point
results and the property (P), Comput. Math. Appl., \textbf{63} (2012),
1361-1368.

[9] G. Gw\'{o}\'{z}d\'{z}-\L ukowska and J. Jachymski,\ IFS on a metric
space with a graph structure and extensions of the Kelisky-Rivlin theorem,
J. Math. Anal. Appl., \textbf{356} (2009), 453-463.

[10] N. Hussain, M. A. Kutbi, S. Khaleghizadeh and P. Salimi, Discussions on
recent results for $\alpha $--$\Psi $-contractive mappings, Abstr. Appl.
Anal., vol. 2014, Article ID 456482, 13 pages, 2014.

[11] M. Klimek and M. Kosek, Generalized iterated function systems,
multifunctions and Cantor sets, Ann. Polon. Math., \textbf{96} (2009), 25-41.

[12] A. Latif, W. Sintunavarat and A. Ninsri, Approximate fixed point
theorems for partial generalized convex contraction mappings in $\alpha $%
-complete metric spaces, Taiwanese J. Math., \textbf{19} (2015), 315-333.

[13] K. Le\'{s}niak, Infinite iterated function systems: a multivalued
approach, Bull. Pol. Acad. Sci. Math., \textbf{52} (2004), 1-8.

[14] E. Llorens-Fuster, A. Petru\c{s}el and J.-C. Yao, Iterated function
systems and well posedness, Chaos Solitons Fractals, \textbf{41} (2009),
1561-1568.

[15] L. M\'{a}t\'{e}, The Hutchinson-Barnsley theory for certain
noncontraction mappings, Period. Math. Hungar., \textbf{27} (1993), 21-33.

[16] M. A. Miandaragh, M. Postolache and S. Rezapour, Approximate fixed
points of generalized convex contractions, Fixed Point Theory Appl., vol.
2013, article 255, 2013.

[17]\ A. Mihail and R. Miculescu, Applications of fixed point theorems in
the theory of generalized IFS, Fixed Point Theory Appl., 2008, Art. ID
312876, 11 pp.

[18] A. Mihail and R. Miculescu, The shift space for an infinite iterated
function system, Math. Rep.(Bucur.), \textbf{61} (2009), 21-31.

[19]\ A. Mihail and R. Miculescu, Generalized IFSs on noncompact spaces,
Fixed Point Theory Appl., 2010, Art. ID 584215, 15 pp.

[20] A. Petru\c{s}el, Iterated function system of locally contractive
operators, Rev. Anal. Num\'{e}r. Th\'{e}or. Approx., \textbf{33} (2004),
215-219.

[21] A. Petru\c{s}el, Fixed point theory and the mathematics of fractals,
Topics in mathematics, computer science and philosophy, 165-172, Presa Univ.
Clujean\u{a}, Cluj-Napoca, 2008.

[22] N. A. Secelean, Iterated function systems consisting of $F$%
-contractions, Fixed Point Theory Appl., 2013, 2013:277.

[23] N. A. Secelean, Countable iterated function systems, Lambert Academic
Publishing, 2013.

[24] N. A. Secelean, Generalized iterated function systems on the space $%
l^{\infty (X)}$, J. Math. Anal. Appl., \textbf{410} (2014), 847-458.

[25] F. Strobin and J. Swaczyna, On a certain generalization of the iterated
function system, Bull. Australian Math. Soc., \textbf{87} (2013), 37-54.

[26] F. Strobin, Attractors of generalized IFSs that are not attractors of
IFSs, J. Math. Anal. Appl., \textbf{422} (2015), 99-108.

\bigskip

University of Bucharest

Faculty of Mathematics and Computer Science

Str. Academiei\ 14, 010014 Bucharest, Romania

E-mail: miculesc@yahoo.com, mihail\_alex@yahoo.com

\end{document}